\documentclass{amsart}
\usepackage{amssymb, amsmath, amscd, bbm, graphicx}
\swapnumbers

\renewcommand{\a}{\alpha}
\renewcommand{\b}{\beta}
\renewcommand{\c}{\gamma}

\providecommand{\e}{\epsilon}

\providecommand{\s}{\sigma}

\newcommand{\R}{\mathbb R}

\newcommand{\Q}{\mathbb Q}

\newcommand{\calC}{\mathcal{C}}

\newcommand{\conj}{\overline}

\newcommand{\laplacian}{\square}

\newcommand{\im}{\operatorname{im}}

\newcommand{\diag}{\operatorname{diag}}
\newcommand{\lk}{\operatorname{lk}}
\newcommand{\rk}{\operatorname{rk}}
\newcommand{\nullity}{\operatorname{nullity}}
\newcommand{\vspan}{\operatorname{span}}

\numberwithin{equation}{section}

\theoremstyle{plain}
\newtheorem{thm}[equation]{Theorem}
\newtheorem{lem}[equation]{Lemma}
\newtheorem{cor}[equation]{Corollary}
\newtheorem{prop}[equation]{Proposition}

\theoremstyle{remark}
\newtheorem{remark}[equation]{Remark}
\newtheorem{ex}[equation]{Example}

\theoremstyle{definition}
\newtheorem{defn}[equation]{Definition}

\begin{document}

\title{Analytic Subdivision Invariants}
\author{Jer-Chin (Luke) Chuang}
\thanks{The author thanks Chris Rasmussen, Brendan Hassett, and especially Robin Forman for many helpful discussions.}

\begin{abstract}
This paper introduces an inner product on chain complexes of finite simplicial complexes that is well-adapted to the harmonic study of subdivisions.  Its definition utilizes a decomposition of the chain spaces that suggests a sequence of subdivision invariants which we show do not all vanish for non-trivial subdivisions.  We exhibit a combinatorial lower bound for these invariants and provide an effective algorithm for their computation.  Unfortunately, these invariants cannot distinguish every subdivision nor do they necessarily increase over successive subdivisions.
\end{abstract}

\maketitle


\section{Introduction}
\begin{figure}
\includegraphics[scale=.4, bb=0 0 578 146]{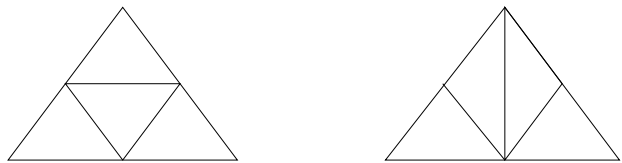}
\caption{Two subdivisions of the 2-simplex with identical face numbers}\label{facevector}
\end{figure}

Let $N$ be a fixed, finite simplicial complex and $M,M'$ subdivisions of $N$.  We say that $M,M'$ are isomorphic provided there is a simplicial isomorphism between them.  A natural question is to determine whether any two given subdivisions are equivalent.  Certainly, one can enumerate all the possible simplicial maps between their vertex schemes, but to the author's knowledge there is not yet a more effective way.  As in topology, one could begin by checking invariants.  By a (real-valued) \emph{subdivision invariant}, we mean an assignment of real numbers $f(M), f(M')$ such that if $M\cong M'$, then $f(M)=f(M')$.  Equivalently, $f$ is a function on the poset of subdivisions of $N$ which descends to the quotient poset of simplicial isomorphism classes.  The most classical examples are the face numbers, but these are easily seen to be insufficient; for example, see Figure \ref{facevector}.

This paper examines a sequence of \emph{analytic} subdivision invariants derived from inner product structures on the chain spaces.  The use of inner product structures to study finite simplicial complexes was initiated by Eckmann in \cite{eckmann45} (see \cite{eckmann00} for more recent applications).  The definition of the invariants is based on a particular decomposition of the chain spaces.  Though each chain space has a canonical inner product given by the simplices, this decomposition can also be used to associate to each subdivision $i\colon N\to M$ an inner product on the chain spaces $C_*(M)$ which reflects the harmonic theory on $C_*(N)$.  The following is an amalgamation of Theorem \ref{decomposition} and Corollary \ref{harmonicrelation}:

\begin{thm}
Let $i\colon N\to M$ be a subdivision of a simplicial complex $N$.  There exists subspaces $V_k\subset C_k(M)$ such that
\begin{equation}
C_k(M) = i_*C_k(N)\oplus V_k \oplus \partial_{k+1}V_{k+1}
\end{equation}
If $h$ is an inner product on $C_*(N)$, then this decomposition defines an inner product $g'$ on $C_*(M)$ which satisfies
\begin{equation}
i_*\laplacian_h = \laplacian_{g'}i_*
\end{equation}
where $\laplacian_h,\laplacian_{g'}$ are the Laplacians relative $h,g'$ respectively.
\end{thm}

This allows the harmonic theory of the standard inner products on $C_*(N),C_*(M)$ to be compared on just one complex.  In particular, one can study the norm of the differences of harmonic representatives.  The connections between harmonic theory and subdivisions will be the subject of a future paper.  Here, the above decomposition is used to define the following sequence of numbers: for each non-negative integer $k$, consider
\begin{equation}
C_k= \max_{|x|=|y|=1} g(x,y)\quad x\in \im i_*, \quad y\in \partial V_{k+1}
\end{equation}
which measures the extent the subspace $\im i_*$ fails to be $\partial^*$-invariant relative the adjoint with respect to the standard inner product on the chain space $C_k(M)$.  The numbers $C_k$ are subdivision invariants (Theorem \ref{invariant}) and do not all vanish for non-trivial subdivisions (Theorem \ref{same}).  The latter claim follows from the following combinatorial bound for these invariants.  To state this bound, we introduce the following definition which will be useful for geometric arguments throughout this paper:

\begin{defn}
Let $i\colon N\to M$ be a subdivision of simplicial complexes.  A \emph{$k$-party} of $M$ is a union of $k$-simplices of $M$ which coincides with the image of a $k$-simplex of $N$ under the simplex map $i$.  The $k$-simplices of $M$ which constitute a $k$-party are called its \emph{members}.
\end{defn}

We sometimes find it expedient to think of parties as linear combinations in the chain space.
\begin{ex}
Consider the $1$-complex which is the union of the two intervals $[0,1]$ and $[1,2]$.  Suppose the former is subdivided at $1/2$.  Then, the resulting complex has two $1$-parties: $[1,2]$ and $\{[0,1/2],[1/2,1]\}$.
\end{ex}

Let $i\colon N\to M$ be a subdivision of a simplicial complex and $\s$ any $k$-simplex of $M$.  Write $F_{\s}$ for the number of incident $(k+1)$-simplices and $N_{\s}$ for the number of $k$-faces of these incident $(k+1)$-simplices which are supported on singly represented parties, i.e. parties that support no two of these $k$-faces.  Here is the combinatorial lower bound (Proposition \ref{lowerbound}):

\begin{prop}
Let $i\colon N\to M$ be a subdivision of a simplicial complex, and suppose that the number of $(k+1)$-simplices has increased.  Pick any $k$-simplex $\s$ not supported on any $k$-parties.  Then,
\begin{equation}
C_k=\langle \im i_*,\partial V_{k+1}\rangle \geq \frac{1}{\sqrt{F_{\s}+(k+1)}}\sqrt{\frac{N_{\s}}{F_{\s}}}
\end{equation}
\end{prop}

These invariants $C_k$ admit a geometric interpretation as a constrained optimization of a quadratic functional over a conic intersection locus.  An effective algorithm for their computation (illustrated in Example \ref{extremalex}) is provided and the codimension one value is deduced for elementary stellar subdivision of an isolated simplex (Proposition \ref{isolatedcase}) or along an interior simplex (Proposition \ref{interiorcase}).  Using the algorithm, one can compute that the sequence $C_k$ distinguishes between the subdivisions in Figure \ref{facevector} whereas the face numbers could not.  However, they are still insufficient to distinguish between any two subdivisions of a common simplicial complex (Remark \ref{insufficient}) and do not necessarily increase with successive subdivision (Example \ref{noincrease}).

The sections are organized as follows: Section \ref{harmonic} introduces an inner product well adapted for harmonic theory on subdivisions via a particular decomposition of the chain spaces.  The sequence of subdivision invariants $C_k$ studied in this paper is introduced in Section \ref{definition} and a combinatorial lower bound is given in Section \ref{properties}.  Section \ref{extremal} provides the geometric interpretation and an effective algorithm for their computation.  This interpretation is used in Section \ref{stellar} to obtain the codimension one $C_k$ values for an elementary stellar subdivision in certain cases.  Finally, the appendix (Section \ref{laplheat}) provides background on the Laplacian for chain complexes.

\section{Subdivisions and Inner Products}\label{harmonic}
Let $i\colon N\to M$ be a subdivision of simplicial complexes.  Suppose $h,g$ are inner products on the chain spaces $C_*(N),C_*(M)$ respectively.  Because the induced map $i_*$ on chain spaces is injective, one can readily define an inner product $g'$ on $C_*(M)$ such that $i^*g'=h$.  Thus, the combinatorial structure of the unsubdivision $N\preceq M$ is encoded at the level of the chain spaces $C_*(M)$.  The choice of $g'$ such that its pullback is $h$ is in general not unique.  In this section, we will show that there is a canonical choice which further satisfies the property that the chain map $i_*$ and the Laplacians commute:
\begin{equation}\label{commute}
i_*\laplacian_h = \laplacian_{g'}i_*
\end{equation}
where $\laplacian_h,\laplacian_{g'}$ are the Laplacians relative $h,g'$ respectively.  This will allow us to conduct Laplacian calculations on the subdivided complex $M$ alone.  Hopefully, questions about a subdivision pair may now be adequately translated into a question about a pair of inner products on a single vector space.  To define the canonical choice, we decompose $C_*(M)$ into convenient summands.

A sufficient, but not necessary condition for Equation (\ref{commute}) is that $i_*$ and $\partial^*_{g'}$ commute.

\begin{prop}\label{suffcommute}
Suppose $i^*g'=h$.  The maps $i_*,\partial_{g'}^*$ commute iff the $g'$-orthogonal complement of $\im i_*$ is $\partial$-invariant.
\end{prop}
\begin{proof}
First, suppose $i_*,\partial_{g'}^*$ commute.  Let $a\in (\im i_*)^{\perp_{g'}}$ and $c\in C_*(N)$ any chain.  Then,
\begin{equation*}
g'(i_*c,\partial a) = g'(\partial_{g'}^* i_*c,a) = g'(i_* \partial_{g'}^* c,a)=0
\end{equation*}
so that $(\im i_*)^{\perp_{g'}}$ is $\partial$-invariant.

Conversely, suppose the $g'$-orthogonal complement of $\im i_*$ is $\partial$-invariant.  For any chain $c\in C_*(M)$ we write $c=i_*b+a$ for $b\in C_*(N)$ and $a\in (\im i_*)^{\perp_{g'}}$.  Then, for any chain $c'\in C_*(N)$,
\begin{equation*}
g'(i_*\partial_{g'}^*c',a)=0=g'(i_*c',\partial a)=g'(\partial_{g'}^* i_*c',a)
\end{equation*}
so that $g'([i_*,\partial_{g'}^*]c',a)=0$.  Also, because $i^*g'=h$ we have
\begin{equation*}
g'(i_*\partial_{g'}^*c',i_*b)=h(c',\partial b)=g'(i_*c',\partial i_*b)=g'(\partial_{g'}^* i_*c',i_*b)
\end{equation*}
so that $g'([i_*,\partial_{g'}^*]c',i_*b)=0$.  Thus, $g'([i_*,\partial_{g'}^*]c',c)=0$ for all chains $c\in C_*(M)$ and $c'\in C_*(N)$ so that $[i_*,\partial_{g'}^*]=0$.
\end{proof}

Note that $\partial$-invariance implies that the induced map
\begin{equation*}
\conj{\partial}_{k+1}\colon \frac{C_{k+1}(M)}{\ker \partial_{k+1}+\im i_*}\to \frac{C_k(M)}{\im i_*}
\end{equation*}
is injective, though this result is true more generally:

\begin{prop}\label{complement}
Let $i\colon N\to M$ be a subdivision of simplicial complexes.  The map
\begin{equation}
\conj{\partial}_{k+1}\colon \frac{C_{k+1}(M)}{\ker \partial_{k+1}+\im i_*}\to \frac{C_k(M)}{\im i_*}
\end{equation}
induced by the boundary map is injective for all $k$.
\end{prop}
\begin{proof}
Consider the diagram
\begin{equation}
\begin{CD}
C_{k+1}(N) @>\partial_{k+1}>> C_k(N)\\
@V i_* VV  @VV i_* V\\
C_{k+1}(M) @>\partial_{k+1}>> C_k(M)
\end{CD}
\end{equation}
Suppose $\partial_{k+1}c\in \im i_*$ for some chain $c\in C_{k+1}(M)$.  Since $\partial^2=0$, the isomorphism $i_*\colon H_k(N)\to H_k(M)$ implies that $i_*^{-1}[\partial c]=0\in H_k(N)$.  Hence, there is a $(k+1)$-chain $a\in C_{k+1}(N)$ such that $i_*^{-1}[\partial c]=[\partial a]$.  We then have $\partial_{k+1}(c-i_*a)\equiv 0$ so that $c\in \ker \partial_{k+1}+\im i_*$.
\end{proof}

\begin{remark}
Note that we do \emph{not} have $(\ker \partial+\im i_*)^{\perp}\to (\im i_*)^{\perp}$.  This is the observation behind the subdivision invariants to be introduced in Section \ref{definition}.
\end{remark}

Now we can write the desired decomposition:
\begin{thm}\label{decomposition}
Let $i\colon N\to M$ be a subdivision of simplicial complexes, and equip $C_*(M)$ with an inner product.  Then, there exist subspaces $V_k\subset C_k(M)$ such that
\begin{equation}\label{decompositioneq}
C_k(M)= (\im i_* \oplus \partial V_{k+1}) \oplus_{\perp} V_k
\end{equation}
\end{thm}
\begin{proof}
We show that $\im i_*+\im \partial_k^M=\im i_*\oplus\partial V_{k+1}$.  Then, defining
\begin{align}\label{canonical}
V_{top} &= (\im i_*)^{\perp_g}\\ 
V_k &= (\im i_* \oplus \partial V_{k+1})^{\perp_g}
\end{align}
yields the desired decomposition.  First, by the homology isomorphism $H_{top}(N)\to H_{top}(M)$, we have $\ker \partial_{top}^M\subset \im i_*$ so that $V_{top}\perp \ker \partial_{top}$ and hence, $V_{top}\perp (\ker \partial_{top}+\im i_*)$.  By Proposition \ref{complement}, $\partial V_{top}\cap \im i_*=0$.  Similarly, the homology isomorphism $H_k(N)\to H_k(M)$ implies that $\ker \partial_k^M\subset \im i_*+\im \partial_{k+1}^M$ so that $V_k\perp (\ker \partial_k+\im i_*)$.  By Proposition \ref{complement} again, $\partial V_{k+1}\cap \im i_*=0$.
\end{proof}

\begin{cor}\label{harmonicrelation}
Let $i\colon N\to M$ be a subdivision of simplicial complexes, and suppose the chain spaces are equipped with inner products $h,g$ respectively.  There exists an inner product $g'$ on $C_*(M)$ such that the associated Laplacians and inclusions commute:
\begin{equation*}
i_*\laplacian_h = \laplacian_{g'}i_*
\end{equation*}
\end{cor}
\begin{proof}
We define an inner product $g'$ on $C_*(M)$ by setting the summands
\begin{equation*}
C_k(M)= \im i_* \perp_{g'} \partial V_{k+1} \perp_{g'} V_k
\end{equation*}
mutually orthogonal while using the induced inner product from $(C_*(M),g)$ on both $\partial V_{k+1},V_k$ but $i^*g'=h$ on $\im i_*$.  By construction $i^*g'=h$ and the $g'$-complement of $\im i_*$ is $\partial$-invariant so that Proposition \ref{suffcommute} implies that $i_*,\partial^*_{g'}$ commute.  Hence, $i_*\laplacian_h = \laplacian_{g'}i_*$.
\end{proof}

We will call the inner product $g'$ defined in the preceding argument the \emph{canonical inner product for the morphism} $i\colon (C_*(N),h)\to (C_*(M),g))$ of inner product spaces.  When we use the standard inner products given by the simplices, we will say that $g'$ is the \emph{canonical inner product for subdivision $i$}.

The following relations summarize the behavior of the decomposition with respect to the boundary map $\partial$ and its adjoint $\partial^*$ relative either inner product $g$ or $g'$:
\begin{alignat}{2}\label{bdymaps}
\partial^*V_{k-1} &= 0 &\quad\quad V_k &\rightleftharpoons \partial V_k\\
\partial(\partial V_{k+1}) &= 0 &\quad\quad \im i_* &\stackrel{\partial}{\to} \im i_*\notag
\end{alignat}

\begin{remark}
The constructions in this section are algebraic and thus valid in the setting of chain complexes of inner product spaces and injective morphisms which induce isomorphisms on homology.
\end{remark}

\section{An Analytic Approach to Subdivision Complexity}\label{definition}
In this section, we introduce a subdivision invariant based on the canonical decomposition (\ref{decompositioneq}) for a subdivision.  The original context was the study of the norm of the difference of cycle representatives of a fixed homology class (see below).

Let $(V,g)$ be an inner product space, and $X,Y$ subspaces.  We write $\langle X,Y\rangle_g$ for the quantity
\begin{equation}
\langle X,Y\rangle_g = \max_{|x|=|y|=1} g(x,y)
\end{equation}
where $x\in X$ and $y\in Y$, and we will omit the subscript when using the standard inner product.

Recall that by Theorem \ref{decomposition}, for each subdivision $i\colon N\to M$, there is an associated canonical decomposition relative the standard inner product on $C_*(M)$,
\begin{equation*}
C_*(M)=(\im i_*\oplus \partial V_{*+1})\oplus_{\perp} V_*
\end{equation*}

\begin{defn}
For each non-negative integer $k$, define
\begin{equation}
C_k=\langle \im i_*,\partial V_{k+1}\rangle
\end{equation}
\end{defn}

\begin{thm}\label{invariant}
The values $C_k$ are invariants of the subdivision.
\end{thm}
\begin{proof}
Let $i'\colon N\to M'$ be a subdivision and $\varphi\colon M\to M'$ a simplicial isomorphism.  Because $\varphi$ is a dimension-preserving map between simplices (in fact, just re-labeling simplices), the induced map $\varphi_*\colon (C_*(M),std_M)\to (C_*(M'),std_{M'})$ is an isometry and we have $V_k'=\varphi_* V_k$.  Thus,
\begin{equation*}
\langle \im i_*',\partial V_{k+1}'\rangle = \langle \varphi_*\im i_*,\varphi_*\partial V_{k+1}\rangle = \langle \im i_*,\partial V_{k+1}\rangle
\end{equation*}
so that the values $C_k$ are identical for isomorphic subdivisions of a given complex.
\end{proof}

\begin{remark}
Evidently, $C_k\in [0,1]$, and if $i$ is the trivial subdivision, then $C_k=0$ for all $k$.  We will see subsequently (Theorem \ref{same}) that these invariants $C_k$ do not all vanish for non-trivial subdivisions.
\end{remark}

These quantities $C_k$ arise somewhat naturally in the study of the norm of the difference between certain cycle representatives of a fixed homology class.  Again, let $i\colon N\to M$ be a subdivision and endow $C_*(M)$ with an inner product $g$.  Fixing a homology class for $|N|$, let $\a$ be the harmonic representative in $C_*(M)$ relative $g$ and $\a'$ the inclusion under $i_*$ of any cycle representative in $C_*(N)$.  Then, from the relations
\begin{equation*}
|\a'-\a|^2_g = g(\a'-\a,\a'-\a)=g(\a',\a'-\a)\leq C|\a'|_g|\a'-\a|_g
\end{equation*}
where $C=\langle i_*\ker \partial^N,\im \partial_*\rangle_g$, we see that
\begin{equation*}
|\a'-\a|_g\leq C |\a'|_g
\end{equation*}
The subspace $\im i_*$ inherits an induced inner product and has a Hodge decomposition:
\begin{equation*}
\im i_* = i_*\ker \partial^N \oplus_{\perp_g} \im \partial_g^*
\end{equation*}
because $\ker \partial\cap \im i_*=i_*\ker \partial^N$ by injectivity of $i_*$.  Then,
\begin{equation*}
C=\langle i_*\ker \partial^N,\im \partial_*\rangle_g=\langle \im i_*, \im \partial_*\rangle_g
\end{equation*}
Since these two subspaces intersect non-trivially, $C=1$ so that
\begin{equation}
|\a'-\a|_g\leq |\a'|_g
\end{equation}
However, from the canonical decomposition, we know that $\im \partial^M_k=i_*\partial^N_k\oplus \partial V_{k+1}$ so that taking $C_k=\langle \im i_*,\partial V_{k+1}\rangle$ recovers a potentially interesting extremization.  In fact, since $\im \partial_k^*\subset \im i_*\oplus_{\perp_g} V_{k+1}$, we know that
\begin{equation*}
C_k=\langle \im i_*,\partial V_{k+1}\rangle=\langle \partial_k^*\im i_*, V_{k+1}\rangle
\end{equation*}
so that $C_k$ is a measure of the failure of $\im i_*$ to be closed under the adjoint $\partial^*$ relative $g$.  Alternatively, it is a measure of how $\im i_*$ and $\partial V_{k+1}$ project onto each other and hence the ``angle'' between them.

\begin{remark}
One could have used instead the canonical inner product $g'$ associated to the subdivision, and similarly obtain $|\a-\a'|_{g'}\leq C'|\a|_{g'}$ where $C'=\langle \ker \partial,\im \partial_{k+1}\rangle_{g'}$, the chain $\a'$ is $g'$-harmonic, and $\a$ is an arbitrary cycle. Again $C'=1$ but because $\im \partial\subset \ker \partial$, there is no interesting extremization in this case.  Indeed, $\im i_*$ is by construction $\partial_{g'}^*$-invariant since the pullback of $g'$ by the subdivision map is the standard inner product on $N$, as if everything were done on $N$ alone.  More easily, this also follows because $i_*,\partial^*_{g'}$ commute.
\end{remark}

\begin{remark}
If in our analysis, we restricted $\a'$ to be $g'$-harmonic, then one can show that the subspaces $\ker \laplacian_{g'}$ and $\im \partial_*$ intersect trivially and that $\gamma = \langle \ker \laplacian_{g'}, \im \partial_*\rangle$ defines another subdivision invariant.  Analysis of this invariant, though only useful in cases of non-trivial homology, will be the subject of a future paper.
\end{remark}

\section{Properties of the Invariants $C_k$}\label{properties}
In this section, we derive a combinatorial lower bound on the invariants $C_k$ (Theorem \ref{lowerbound}) which will show that these invariants do not all vanish for non-trivial subdivisions (Theorem \ref{same}).  We also discuss the difficulties in estimating $C_k$ over successive subdivisions and provide an example where $C_0$ \emph{does not} increase with successive subdivision (Example \ref{noincrease}).

\subsection{Combinatorial Lower Bound for Invariants $C_k$}
The central result of this subsection is a combinatorial lower bound for the invariants $C_k$ (Theorem \ref{lowerbound}).  This bound implies that for any non-trivial subdivision, there is a $C_k\in (0,1)$ (Theorem \ref{same}).  We also present explicit bases for the subspaces $V_1,V_2$ (Propositions \ref{isbasis1} and \ref{isbasis2}).  We begin by computing the dimension of $V_k$.

\begin{prop}\label{dimV}
Let $s_k$ denote the number of $k$-simplices of a simplicial complex.  If $i\colon N\to M$ is a subdivision, then
\begin{align}
\dim V_0 &= 0\\
\dim V_k &= (-1)^k \left(\sum_{j=0}^{k-1} (-1)^j s_j(N) - \sum_{j=0}^{k-1} (-1)^j s_j(M)\right)\quad\quad k\geq 1
\end{align}
In particular, $\dim V_1=s_0(M)-s_0(N)$ the number of new vertices.
\end{prop}
\begin{proof}
From Theorem \ref{decomposition}, we have $C_k(M)=\im i_*\oplus \partial V_{k+1}\oplus V_k$.  Then,
\begin{equation*}
s_k(M)=s_k(N)+\rk \partial_{k+1}^V + \dim V_k
\end{equation*}
where $\partial_{k+1}^V$ is the restiction of $\partial_{k+1}$ to $V_{k+1}$.  Now $\dim V_k=\rk \partial_k^V + \nullity(\partial_k^V)$.  By Proposition \ref{complement}, we have $\nullity(\partial_k^V)=0$ so that $\dim V_k = \rk \partial_k^V$ and
\begin{equation*}
s_k(M)= s_k(N) + \rk \partial_{k+1}^V + \rk \partial_k^V
\end{equation*}
Then, taking an alternating sum telescopes the $\rk \partial_k^V$ terms:
\begin{equation*}
\sum_{j\geq k} (-1)^j s_j(M) = \sum_{j\geq k} (-1)^j s_j(N) + (-1)^k \rk \partial_k^V = \sum_{j\geq k} (-1)^j s_j(N) + (-1)^k \dim V_k
\end{equation*}
Since
\begin{equation*}
\sum_{j\geq 0} (-1)^j s_j(M)= \chi(M)=\chi(N)=\sum_{j\geq 0} (-1)^j s_j(N)
\end{equation*}
we have
\begin{align*}
\dim V_0 &= 0\\
\dim V_k &= (-1)^k \left(\sum_{j=0}^{k-1} (-1)^j s_j(N) - \sum_{j=0}^{k-1} (-1)^j s_j(M)\right)\quad\quad k\geq 1
\end{align*}
as claimed.
\end{proof}

\begin{remark}
Since $\dim V_k\geq 0$, the preceding proposition implies the inequality $\sum_{j=0}^{k-1} (-1)^j s_j(M)\geq \sum_{j=0}^{k-1} (-1)^j s_j(N)$ for $k$ odd and the reverse inequality for $k$ even.
\end{remark}

\begin{lem}\label{inV}
If $\s$ is a $k$-simplex not supported on any $k$-party, then $\partial^*\s\in V_{k+1}$.
\end{lem}
\begin{proof}
By the relations in (\ref{bdymaps}), $\im \partial^*\subset \im i_*\oplus_{\perp}V_{k+1}$, so we need only show that $\partial^*\s$ is perpendicular to $\im i_*$.  We present two arguments, the first analytic, the second geometric.

\emph{Analytic argument}: For any $a\in C_{k+1}(N)$, we have $(\partial^*\s,i_*a)=(\s,i_*\partial a)=0$ by hypothesis.  Thus, $\partial^*\s$ is perpendicular to $\im i_*$ and hence an element of $V_{k+1}$.

\emph{Geometric argument}: First, the $(k+1)$-chain $\partial^*\s$ is supported on $(k+1)$-simplices incident to $\s$. Now, either the $k$-simplex $\s$ is supported on a $(k+1)$-party or not.  If not, then $\partial^*\s$ is not supported on any $(k+1)$-party and hence perpendicular to $\im i_*$.  If so, incident $(k+1)$-simplices are of two types depending on whether they are members of the necessarily non-singular $(k+1)$-party.  Now, the party members are $(k+1)$-dimensional simplices of a subdivision of a $(k+1)$-simplex.  Hence, there are only two such, namely those determined by having $\s$ as their shared boundary, and their contribution in $\partial^*\s$ is thus perpendicular to $\im i_*$.  The remaining incident $(k+1)$-simplices are not supported on $\im i_*$, and we conclude that $\partial^*\s$ is perpendicular to $\im i_*$.
\end{proof}

Let $i\colon N\to M$ be a subdivision of a simplicial complex and $\s$ any $k$-simplex of $M$.  Write $F_{\s}$ for the number of incident $(k+1)$-simplices and $N_{\s}$ for the number of $k$-faces of these incident $(k+1)$-simplices which are supported on singly represented parties, i.e. parties that support no two of these $k$-faces.

\begin{prop}\label{lowerbound}
Let $i\colon N\to M$ be a subdivision of a simplicial complex and suppose that the number of $(k+1)$-simplices has increased.  Pick any $k$-simplex $\s$ not supported on any $k$-parties.  Then,
\begin{equation}\label{bound}
C_k=\langle \im i_*,\partial V_{k+1}\rangle \geq \frac{1}{\sqrt{F_{\s}+(k+1)}}\sqrt{\frac{N_{\s}}{F_{\s}}}
\end{equation}
\end{prop}
\begin{proof}
Since the number of $(k+1)$-simplices has increased, there is a non-singular $(k+1)$-party.  Any two party members will share a $k$-simplex $\s$ not supported on any $k$-party.  By the Lemma \ref{inV}, we have $\partial^*\s\in V_{k+1}$ so that extremization over the subspace $\langle \partial\partial^*\s\rangle\subset \partial V_{k+1}$ gives a lower bound
\begin{equation*}
C_k=\langle \im i_*,\partial V_{k+1}\rangle \geq \langle \im i_*,\langle\partial\partial^*\s\rangle\rangle
\end{equation*}
Now, $\partial\partial^*\s=F_{\s} \s+\sum \b$ where the latter sum has $\mu=((k+2)-1)F_{\s}=(k+1)F_{\s}$ terms,  the number of $k$-faces of $\partial^*\s$ other than $\s$.  Then, $|\partial\partial^*\s|^2=F_{\s}^2+F_{\s}(k+1)$ and $|\pi\partial\partial^*\s|^2=N_{\s}$ where $\pi$ is the projection onto $\im i_*$.  The latter equality reflects the fact that if two members of the same $k$-party are present in $\partial\partial^*\s$ (and by the geometry of simplices, there can be at most two), then their contribution is orthogonal to $\im i_*$.  Thus,
\begin{equation*}
C_k^2\geq \frac{|\pi(\partial\partial^*\s)|^2}{|\partial\partial^*\s|^2} = \frac{N_{\s}}{F_{\s}^2+F_{\s}(k+1)}
\end{equation*}
and the inequality follows.
\end{proof}

\begin{ex}\label{C0stellarcase}
Let $i\colon N\to M$ be an elementary stellar subdivision with $v$ being the new vertex.  Write $N$ for the number of vertices in its link $\lk(v,M)$.  Then, $F_v=N=N_v$ and since $\dim V_1=1$, we have $\langle\partial\partial^*v\rangle=V_1$ so that
\begin{equation}
C_0=\langle \im i_*,\partial V_1\rangle = \frac{1}{\sqrt{N+1}}
\end{equation}
In particular, all elementary stellar subdivisions of an isolated $d$-simplex have the same $C_0$.  This equality offers a partial combinatorial description of the analytic invariant $C_0$.  The invariant $C_0$ also gives a lower bound on $N_v$, but this is useful only when $C_0>1/2$.
\end{ex}

\begin{ex}
For the elementary stellar subdivision in Figure \ref{stellarex}, choosing the edge $e_5$, we have $F=2$ and $N=2$ (because the contribution in $\partial\partial^*e_5$ by the edges in the base $1$-party will be perpendicular to $\im i_*$) so that
\begin{equation}
C_1=\langle \im i_*,\partial V_2\rangle \geq \frac{1}{\sqrt{2+(1+1)}}=\frac{1}{2}
\end{equation}
and $C_0=1/2$ by the preceding example.
\end{ex}

\begin{ex}\label{combmfdcase}
Consider the codimension one invariant for a subdivision of a combinatorial $d$-manifold.  For any $(d-1)$-simplex $\s$, we have $F_{\s}=2$.  Hence,
\begin{equation*}
C_{d-1}=\langle \im i_*,\partial V_d\rangle \geq \max \frac{1}{\sqrt{2+d}}\sqrt{\frac{N_{\s}}{2}}=\sqrt{\frac{\max N_{\s}}{2(d+2)}}
\end{equation*}
where the maximum is taken over all codimension one simplices $\s$.  When $C_{d-1}$ is known, this gives a lower bound for $N_{\s}$.
\end{ex}

\begin{thm}\label{same}
If the number of $(k+1)$-simplices has changed due to a subdivision, then $C_k=\langle \im i_*,\partial V_{k+1}\rangle\in (0,1)$.  In particular, for non-trivial subdivisions of pure simplicial complexes, we have $C_k\in (0,1)$ in all positive codimensions.
\end{thm}
\begin{proof}
This follows from the positivity of the lower bound under the hypothesis.  Subdivisions of pure simplicial complexes increase the number of simplices in any dimension where simplices were already present.
\end{proof}

In the remainder of this subsection, we present explicit descriptions for the subspaces $V_1,V_2$.

\begin{prop}\label{isbasis1}
Let $\{v_j\}$ be the set of new vertices.  Then, $\{\partial^*v_j\}$ is a basis for $V_1$.
\end{prop}
\begin{proof}
Since the new vertices $v_j$ are not supported on $0$-parties, by Lemma \ref{inV}, we have $\partial^*v_j\in V_1$.

Next, we show that the $1$-chains $\partial^* v_j$ are linearly independent.  Note that we need only show that $\im \partial_1$ is not orthogonal to any non-zero subspace of the span of the new vertices.  Let $\s=\sum a_j v_j$ and suppose $a_k\neq 0$.  Let $\tau$ be any $1$-chain that is a directed path from some original vertex to $v_k$.  Then, $(\partial\tau, \s)=(v_k,v_k)\neq 0$ so that the $1$-chains $\partial^* v_j$ are linearly independent.  Since by Proposition \ref{dimV},  $\dim V_1$ is the number of new vertices, we conclude that $\{\partial^*v_j\}$ is a basis for $V_1$.
\end{proof}

Now, we describe $V_2$ by exhibiting a basis $\{\partial^*e_j\}$ where $e_j$ are particular $1$-simplices.  Note that
\begin{equation*}
\dim V_2=(s_0(N)-s_1(N))-(s_0(M)-s_1(M))=(s_1(M)-s_1(N))-\Delta_0
\end{equation*}
where $\Delta_0=s_0(M)-s_0(N)$ is the number of new vertices.  We will specify $s_1(M)-\dim V_2=\Delta_0+s_1(N)$ number of $1$-simplices, and the remaining edges will then be our desired $e_j$.  First, write $\Delta_0=T+F$ where $T$ is the number of new vertices supported on $1$-parties and $F$ the number not supported.  Then, the number of $1$-simplices supported on $1$-parties is $s_1(N)+T$.  For each $k$-party with $k\geq 2$, choose a spanning tree for the (new) interior vertices, and pick any edge supported on that party connecting the tree to an old vertex, which we will call the \emph{anchor} for the tree.  The number of chosen edges in each $k$-party is then the number of new vertices supported in its interior, and hence $F$ edges have been specified by the trees and anchors.  Together with edges in $1$-parties, we have now specified $(s_1(N)+T)+F=s_1(M)-\dim V_2$ edges.  The remaining $(\dim V_2)$ number of edges are our $e_j$.

\begin{prop}\label{isbasis2}
Let $e_j$ be edges as described above.  Then, $\{\partial^*e_j\}$ is a basis for $V_2$.
\end{prop}
\begin{proof}
By construction, such edges are not supported on $1$-parties.  Hence, we conclude by Lemma \ref{inV} that $\partial^*e_j\in V_2$.

To show that the set $\{\partial^*e_j\}$ is linearly independent, we proceed as we did for $V_1$, namely we will show that no non-trivial subspace of the span of $\{e_j\}$ is orthogonal to $\im \partial_2$.  first, note that for any such edge $e_j$, either both endpoints are on the same spanning tree or not.

If so, let $S$ denote the set of edges in the tree as well as the edge connection to the anchor.  Then $e_j$ determines a $1$-cycle in conjunction with the spanning tree.  This cycle is supported on the $k$-party defining the tree, and since $k$-parties are topologically trivial, the cycle is a boundary.  Hence, there is a $2$-chain $f_j$ supported on the $k$-party for which $\partial f\equiv e_j \mod \langle S\rangle$.

If the endpoints are on different spanning trees, let $R,S$ denote there respective edges and anchor connections, and $P$ the party of minimal dimension supporting both $R,S$.  (The party $P$ exists because $e_j$ is not supported on any $1$-parties.)  On each spanning tree, there is a (unique) path from the supported vertex of $e_j$ to the respective anchor.  The anchors themselves are connected by $1$-parties because they are old vertices, and this determines a $1$-cycle supported on $P$.  Again, topological triviality implies the existence of a supported $2$-chain $f_j$ such that now $\partial f_j\equiv e_j \mod \langle R,S,\im i_*\rangle$.  Hence for each edge $e_j$, there is a $2$-chain $f_j$ such that $(\partial f_j,e_j)=(e_j,e_j)\neq 0$ so that no non-trivial subspace of $\{e_j\}$ is orthogonal to $\im \partial_2$.
\end{proof}

\subsection{Successive Subdivisions}
In this subsection we describe the difficulty in estimating the invariants $C_k$ over successive subdivisions and provide an example (\ref{noincrease}) that shows $C_0$ need not increase with successive subdivision.

Let $X\stackrel{i}{\to} Y\stackrel{j}{\to} Z$ be a sequence of subdivisions.  The fundamental observation is that the canonical inner product for the entire subdivision differs from that obtained via sequential application.  With respect to the subdivision $j\circ i$ we have
\begin{equation}\label{simul}
C_*(Z) = (\im (j\circ i)_*\oplus \partial V_{k+1}) \perp V_k
\end{equation}
whereas for $j,i$ individually, we have
\begin{equation}\label{successive}
C_*(Z) = (\im j_*\oplus \partial B_{k+1})\perp B_k \quad\quad C_*(Y) = (\im i_* \oplus \partial A_{k+1})\perp A_k
\end{equation}
where all the orthogonality relations are with respect to the standard inner product on the appropriate ambient space.  The crux of the discrepancy is that inclusion by $j_*$ of $C_*(Y)$ into $C_*(Z)$ does not preserve the orthogonality relation originally present in $C_*(Y)$.  Hence, we only have
\begin{equation}\label{succ-decomp}
C_*(Z) = (\im (j\circ i)_*\oplus \partial(j_*A_{k+1}\oplus B_{k+1})\oplus j_*A_k) \perp B_k
\end{equation}
In particular, $j_*A_k$ may not be orthogonal to the subspace $\im (j\circ i)_*$ because the complements $A_*$ are determined by the inner product on the intermediate $Y$, not by the standard inner product on $C_*(Z)$.  Moreover, neither do we necessarily have $\partial j_*A_{k+1}\subset \partial V_{k+1}$.  To obtain a decomposition respecting the $V_*$ subspaces, we would need to use $(Y,j^*g)$ the pullback of the standard inner product on $C_*(Z)$.  Then, $V_k = j_*A_k\perp B_k$.

\begin{remark}
Let $a\in A_{k+1}$ and $b\in B_{k+1}$.  Then, $(\partial_*j_*a,\partial b)=(j_*a,\partial^*\partial b)=0$ because $\partial^*\partial B\subset B_{k+1} \perp \im j_*$.  Thus, we may refine the decomposition (\ref{succ-decomp}) as
\begin{equation}
C_*(Z) = (\im (j\circ i)_*\oplus (\partial j_*A_{k+1}\perp \partial B_{k+1})\oplus j_*A_k) \perp B_k
\end{equation}
\end{remark}

\begin{ex}
Let $X=[0,1]$, $Y$ be a subdivision by adding a vertex at $1/2$, and $Z$ another vertex at $1/4$, which we summarize as follows:
\begin{equation*}
X=[0,1]\stackrel{i}{\to} Y=[0,1/2,1] \stackrel{j}{\to} Z=[0,1/4,1/2,1]
\end{equation*}
Numbering the edges in increasing order within each complex, we obtain the following relations in the top dimension:
\begin{alignat*}{2}
\im i_* &= \langle (1,1)\rangle  & A_1 &= \langle (-1,1)\rangle\\
\im j_* &= \langle (1,1,0), (0,0,1)\rangle & B_1 &= \langle (-1,-1,0)\rangle\\
\im (j\circ i)_* &= \langle (1,1,1)\rangle & V_1 &= \langle (-1,1,0),(-1,0,1)\rangle\\
\end{alignat*}
Thus,
\begin{equation*}
j_*A_1 =
\vspan \left\{\begin{pmatrix}1 & 0\\ 1 & 0\\ 0 &1\end{pmatrix}
\begin{pmatrix}-1\\ 1\end{pmatrix}\right \}
=\vspan \left \{\begin{pmatrix}-1\\ -1\\ 1\end{pmatrix}\right \}
\end{equation*}
but $(-1,-1,1)\cdot (1,1,1)\neq 0$ so that $j_*A_1$ is not orthogonal to $\im (j\circ i)_*$.

This example also shows that $j_*\partial A_1$ may not be a subset of $\partial V_1$.  If we label the vertices in increasing order, then
\begin{equation*}
\partial j_*A_1=\langle (1,0,-2,1)\rangle \quad\quad \partial V_1=\langle (1,-2,1,0),(1,-1,-1,1)\rangle
\end{equation*}
and one readily checks that $(1,0,-2,1)\notin \partial V_1$.
\end{ex}

Using the notation of Equations (\ref{simul}) and (\ref{successive}), the question whether the invariants increase with successive subdivision is then a comparison of the quantitites:
\begin{equation}
\langle \im i_*,\partial A_{k+1}\rangle_Y \quad\quad \langle \im (j\circ i)_*,\partial V_{k+1}\rangle_Z
\end{equation}
where the subscripts $Y,Z$ refer to the standard inner products on $C_*(Y), C_*(Z)$ respectively.  The following example shows that not all the invariants necessarily increase with successive subdivision.

\begin{ex}\label{noincrease}
We show that, unfortunately, the invariant $C_0$ does not necessarily increase over successive subdivisions.  Consider a simplicial complex under two successive elementary stellar subdivisions.  Let $v,w$ be the added vertices in that order.  Suppose that their star neighborhoods are disjoint in the final complex so that all the vertices in their links are old vertices.  Let $F_v,F_w$ be the number of vertices in their respective links.  Then, $\partial\partial^*v=(F_v v-\s_v)$ where $\s_v=\sum_{x\in \lk(v)} x$ and likewise for $\partial\partial^*w$.  By Proposition \ref{isbasis1}, $\pi\partial V_1=\langle -\s_v, -\s_w\rangle$ where $\pi$ is the projection onto the inclusion of original chains.  Set
\begin{equation*}
\s = a(F_v v-\s_v) + b(F_w w-\s_w)
\end{equation*}
Then, $C_0$ for the combined subdivision is the maximum of $\pi(\s/|\s|)$ over all real $a,b$, which we seek by equivalently maximizing its square, the function
\begin{equation*}
f(a,b) = \frac{a^2F_v+b^2F_w}{a^2(F_v^2+F_v)+b^2(F_w^2+F_w)}
\end{equation*}
where $(a,b)\neq (0,0)$.  Note that the function only depends on the ratio $b/a$.  Taking partials, we find that
\begin{equation*}
\frac{\partial f}{\partial a} = \frac{ab^2(F_v-F_w)}{(a^2(F_v^2+F_v)+b^2(F_w^2+F_w))^2}
\quad\quad
\frac{\partial f}{\partial b} = \frac{ba^2(F_w-F_v)}{(a^2(F_v^2+F_v)+b^2(F_w^2+F_w))^2}
\end{equation*}
If the links have the same number of vertices, $F_v=F_w$, then the partials vanish, and we find that $f(a,b)=1/(N+1)$ except at $(0,0)$.  Hence, $C_0=1/\sqrt{N+1}$ the same as a single elementary stellar subdivision.
\end{ex}

\section{Computing the Invariants $C_k$}\label{extremal}

In this section, we redescribe the invariants $C_k$ as the maximum of certain quadratic functionals over the (conic) intersection locus of an ellipsoid and hyperplanes.  This interpretation provides another argument (Theorem \ref{anothersame}) for the non-triviality of the invariants $C_k$ earlier shown (Theorem \ref{same}).  The constrained optimization is amenable to the method of Lagrange multipliers, and we show that the desired extremal value corresponds to one of the multipliers (Proposition \ref{simplification}).  We then provide an effective algorithm for its computation which uses only linear algebra over a single-variable polynomial ring.  This is illustrated in Example \ref{extremalex}.

\subsection{Geometric Interpretation of the Invariants $C_k$}\label{geometric}
By the definition of the standard inner product, for $x\in \im i_*$ and $y\in \partial V_{k+1}$, we have
\begin{equation}
g(x,y) = g(x,Lz)
\end{equation}
where $L=\pi\circ\partial\colon V_{k+1}\to \im (i_k)_*$ is the composition of $\partial$ and the ``projection'' $\pi$ onto the unsubdivided $k$-chain subspace $i_*C_k(N)$.  (Note that $\ker \partial\subset \ker L$.)  Explicitly, let $\{\b_j\}$ be the $(k+1)$-simplices of the subdivided complex, and $\{\gamma_i\}$ the $k$-parties.  Then, $L=(L_i^j)$ is the restriction to $V_k$ of the linear map
\begin{equation}
L_i^j = (\gamma_i,\partial \b_j)=
\begin{cases}
\pm 1 \quad \text{if} \quad \partial\b_j \quad \text{supports} \quad \gamma_i \quad\\
0 \quad \text{else}
\end{cases}
\end{equation}
where the sign is given by the relative orientations of $\gamma_i$ and its term in $\partial\b_j$.  Choosing $x\in \im i_*$ in the direction $Lz$ that maximizes $|Lz|$, we see that
\begin{equation}
\langle \im i_*,\partial V_{k+1}\rangle = \max_{|x|=|y|=1} g(x,y) = \max_{|\partial z|=1} |Lz|
\end{equation}
where $z\in V_{k+1}$.  Equivalently, we are computing the maximum singular value of $L$.  The advantage of projecting onto $\im i_*$ rather than $\partial V_{k+1}$ is that there is a natural basis for $\im i_*$, and it is easily written with respect to the standard basis of the subdivided simplices.

Consider first the codimension one case.  Recall that $C_{top}(M)=\im i_*\perp_g V_{top}$.  The condition $z\in V_{top}$ is equivalent to a set of hyperplane conditions given by the top-dimensional parties.  For example if $\{\b_1,\ldots,\b_p\}$ constitute a party and $b_1,\ldots,b_p$ are the coordinate variables for these simplices, then the associated hyperplane condition is $\e_1 b_1+\cdots+\e_p b_p=0$, where $\e_i=\pm 1$.  Next, the condition $|\partial z|=1$ is an ellipsoid (because $\partial$ is not necessarily orthogonal) so that this extremization may be geometrically interpreted as maximizing the norm $|Lz|$ on the intersection locus of an ellipsoid and a set of hyperplanes.  Note that because the boundary map $\partial$ may have non-trivial kernel, the ellipsoid is generally \emph{not} codimension one.  In coordinates, writing $\{\b_i\}$ for the top-dimensional simplices and $P_j$ for the party hyperplanes, we are maximizing $|Lz|$ over the variety (a conic)
\begin{equation}
\mathbf{V}(|\partial z|^2-1, P_j)
\end{equation}
where $z=\sum b_i\b_i$.

For higher codimension, say dimension $k$, since
\begin{equation*}
C_k(M) = (\im i_* \oplus \partial V_{k+1}) \oplus_{\perp} V_k
\end{equation*}
we would use $(k+1)$-simplices $\b_i$ as a basis and obtain hyperplanes $P_j$ due to the party-relations $V_k\perp \im i_*$ and also hyperplanes $Q_l$ due to the relations $V_{k+1}\perp \partial V_{k+2}$.  The latter could be defined by choosing a basis $\{\gamma_l\}$ for $\partial V_{k+2}$ in terms of $\{\b_i\}$ and setting $Q_l$ to be the orthogonality relations $V_{k+1}\perp \langle\gamma_l\rangle$.  However, the nature of the extremization problem remains identical.  (See Example \ref{extremalex}.)

Using the preceding formulation, we now present another argument for the non-triviality of the invariants $C_k$ which we first demonstrated in Theorem \ref{same} and restate here for convenience.

\begin{thm}[Restatement of Theorem \ref{same}]\label{anothersame}
If the number of $(k+1)$-simplices has changed due to a subdivision, then $C_k=\langle \im i_*,\partial V_{k+1}\rangle\in (0,1)$.  In particular, for non-trivial subdivisions of pure simplicial complexes, we have $C_k\in (0,1)$ in all positive codimensions.
\end{thm}
\begin{proof}
Let $k$ be of positive codimension, and $P$ be a non-singular $(k+1)$-dimensional party (the existence of which is assured by hypothesis).  Let $\b_1$ be any member of $P$ incident to $\partial P$.  Since $\b_1$ cannot be incident to all $k$-dimensional parties of $\partial P$, let $\b_2$ be a member of $P$ incident to any yet unrepresented $k$-dimensional party of $\partial P$.  Set $\s=\e_1\b_1+\e_2\b_2$ where $\e_1$ is the coefficient for $\b_1$ in the party definition of $P$ but $\e_2$ is the negative of that for $\b_2$.  By construction $\s\cdot P=0$ so that $\s\perp \im i_*$, and since there is some $k$-dimensional party of $\partial P$ incident to precisely one of $\b_1,\b_2$ we have $L\s\neq 0$.  Now, \emph{a priori} $\s$ is only normal to $\im i_*$ and may not be in $V_{k+1}$.  But we can modify $\s$ by an element of $\partial V_{k+2}$ so that this is true.  This modification does not affect $L\s$ since the definition of $L$ involves composition with $\partial$.  In particular, $0\notin \s+\partial V_{k+2}$.
\end{proof}

\begin{figure}
\includegraphics[scale=.5, bb=47.6 483.42025 542.11428 635.82857]{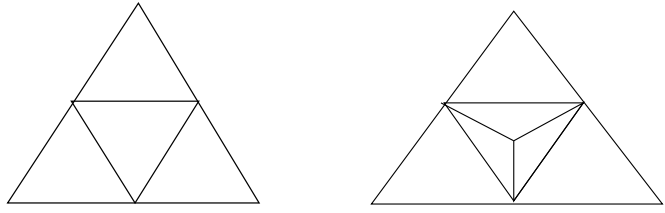}
\caption{Subdivisions of the 2-simplex distinguished by the invariants $C_k$}\label{intstellar}
\end{figure}

\begin{remark}
The extremal values under study can still detect subdivisions away from $\im i_*$, because though the functional $|Lz|^2$ has not changed, the intersection locus has.  Figure \ref{intstellar} provides an example for two subdivisions of the $2$-simplex.  The less subdivided has a value of $1/3$ while the more subdivided a value of $4/11$.
\end{remark}

Let $f=|Lz|^2$, $g=0$ the ellipsoid, and $h_j=0$ the hyperplanes.  When the vectors $\nabla g,\nabla h_j$ are linearly independent over the intersection locus $g=h_j=0$, the method of Lagrange multipliers provides a necessary condition for the critical points of this constrained optimization, namely
\begin{align}\label{system}
\nabla f &= \lambda \nabla g + \sum \lambda_j \nabla h_j\\ \notag
g &= 0\\
h_j &= 0\notag
\end{align}

We will find the following classical result useful:
\begin{quote}
\emph{Euler's Identity for Homogeneous Functions}:  Let $P$ be a homogeneous polynomial over $x_i$, and set $x=(x_1,\ldots,x_n)$.  Then,  $x\cdot \nabla P = (\deg P) P$.
\end{quote}

\begin{prop}\label{linind}
The vectors $\nabla g,\nabla h_j$ are linearly independent over the intersection locus.
\end{prop}
\begin{proof}
Recall that the hyperplanes $h_j=0$ are of two classes: those denoted $P_j$ due to party relations $V_{k+1}\perp \im i_*$ and those denoted $Q_j$ due to $V_{k+1}\perp \partial V_{k+2}$.  The vectors $\nabla P_i$ are linearly independent because they are precisely the party definitions.  Note that they are elements of $\im i_*$.  Also, since $Q_j$ is given by $V_{k+1}\perp \gamma_j$ for basis $\{\gamma_j\}$ of $\partial V_{k+2}$, the vectors $\nabla Q_j$ are precisely $\gamma_j$ and hence linearly independent.  Next, because $\im i_*$ and $\partial V_{k+2}$ are direct summands, we know that the vectors $\nabla h_j$ (being the combined collection of $\nabla P_i$ and $\nabla Q_j$) are linearly independent.

Finally, we use Euler's identity to show that $\nabla g,\nabla h_j$ are linearly independent over the intersection locus:  Let $D=\nabla g-\sum a_j\nabla h_j$ for any real numbers $a_j$, and set $\tilde{g}=g+1$, the homogeneous part of the ellipsoid condition $g=0$.  If $z$ is a point in the intersection locus, then by Euler's Identity,
\begin{align*}
z\cdot D
&= z\cdot \nabla g - \sum a_j z\cdot\nabla h_j=2\tilde{g}(z)-\sum a_j h_j(z)\\
&= 2(g(z)+1)=2
\end{align*}
Thus, $\nabla g, \nabla h_j$ cannot be linearly dependent anywhere on the intersection locus.
\end{proof}

The following result simplifies our extremization task.  Its proof was suggested to the author by Brendan Hassett.

\begin{prop}\label{simplification}
Let $I$ be the ideal in $\Q[\b_i,\lambda_j,\lambda]$ defined by the system $(\ref{system})$.  Then, $|Lz|^2 \equiv \lambda \mod I$.  Thus, the square of the invariant is the maximum of $\lambda$ on $\mathbf{V}(I)$.
\end{prop}
\begin{proof}
We again invoke Euler's Identity.  Write $z=(\b_i)$ and $g=\tilde{g}-1$ where $\tilde{g}$ is now a homogeneous polynomial.  From the Lagrange multiplier system,
\begin{equation*}
\nabla f = \lambda \nabla g + \sum \lambda_j h_j
\end{equation*}
so that
\begin{align*}
z\cdot \nabla f
&= \lambda z\cdot \nabla \tilde{g} + \sum \lambda_j z\cdot \nabla h_j\\
&= \lambda (\deg \tilde{g}) \tilde{g} + \sum \lambda_j (\deg h_j) h_j\\
&= \lambda (\deg \tilde{g}) (g+1) + \sum \lambda_j (\deg h_j) h_j\\
&\equiv (\deg \tilde{g}) \lambda\quad \mod I
\end{align*}
Since $\deg f=\deg \tilde{g}$, we are done.
\end{proof}

\subsection{Algorithm for Computing the Invariants $C_k$}
The preceding subsection showed that the (squares of the) invariants $C_k$ are the maximum values of certain quadratic functionals over conic intersection loci, a computation amenable to the method of Lagrange multipliers.  In this subsection we present an effective algorithm for their computation using only linear algebra over a single-variable polynomial ring.

There exist many algorithms for solving polynomial systems such as that in the Lagrange multiplier system of Equation (\ref{system}).  For example, one could use elimination theory by computing a Grobner basis over $\Q[\lambda_j,\b_i,\lambda]$ for the ideal $I$ defined by the system $(\ref{system})$ under the lexicographic ordering $\lambda_j\succ \b_i\succ \lambda$.  Since $I\cap \Q[\lambda]$ is principal, by the preceding proposition (\ref{simplification}), the roots of the generator for this elimination ideal are candidate extremal values.  We then try lifting each root to the variety $\mathbf{V}(I)$, and the largest that lifts is the desired value.  Note that if we knew that all candidate extremal values always lift, then the desired extremal value is precisely the largest root of the principal generator.

Alternatively, if we know \emph{a priori} that the preimage of the extremal values is a finite set, then there are robust eigenvalue methods for calculating the $\lambda$ values that are actually attained.  This finiteness condition is not always true, e.g. the elementary stellar subdivision of an isolated simplex (Proposition \ref{isolatedcase}), but it is unclear whether this finiteness may be generic in some appropriate sense.

The following is a linear algebraic algorithm valid in any positive codimension for finding the maximum value of $f=|Lz|^2$.  It is based on the observation that except for the constraint $g=0$, the system $(\ref{system})$ is linear over $\b_i, \lambda_j$ with coefficients at most linear in $\lambda$.

\emph{Algorithm for the maximum values of $|Lz|^2$}:
\begin{enumerate}
\item We re-write the vector equation
\begin{equation}
\nabla f - \lambda\nabla g = \sum \lambda_j\nabla h_j
\end{equation}
from the Lagrange multiplier system (\ref{system}) by defining a vector $\mu=(\lambda_j)$ and a matrix $A$ with columns $\nabla h_j$.  Then, the preceding equation reads
\begin{equation*}
A\mu = \nabla f - \lambda\nabla g
\end{equation*}
By Proposition \ref{linind} the vectors $\nabla h_j$ are linearly independent so that the constant matrix $A$ has a left-inverse.  Thus, the multipliers $\lambda_j$ are expressible in terms of $\lambda$ and $\b_i$, and we may eliminate them.
\item Leaving out the ellipsoid constraint $g=0$, we have a system over $\b_i$,
\begin{equation}
(M_{\lambda})z = 0
\end{equation}
where the entries of $M_{\lambda}$ are linear in $\lambda$.
\item The ellipsoid constraint $g=0$ implies that any $\lambda$ satisfying the Lagrange multiplier system $(\ref{system})$ must cause $M_{\lambda}$ to have non-trivial kernel since a trivial partial solution $(\b_i)=0$ does not lie on the ellipsoid.  Hence,
\begin{equation}
Q(\lambda) := \det M_{\lambda} = 0
\end{equation}
For each real root of $Q(\lambda)$, we may compute the associated kernel of $M_{\lambda}$.
\item A partial solution $\lambda$ lifts to a partial solution $(\lambda,\b_i)$ iff $\ker M_{\lambda}$ intersects the ellipsoid.  This can be determined by projecting any basis for the kernel onto the axial directions.  These are given by the eigenspaces of non-zero eigenvalues of the symmetric operator $\partial^T_{top}\partial_{top}$.
\item Finally, the left-inverse in Step (1) implies that any partial solution $(\lambda,\b_i)$ lifts to a solution of the Lagrange multiplier system (\ref{system}).
\end{enumerate}

\begin{remark}\label{alternative}
In the codimension one case, we could do the following instead of Step (1) above: For each party $j$, select a member $\b_i^{(j)}$ and use the equations
\begin{equation}\label{eliminate}
\frac{\partial f}{\partial \b^{(j)}_i} = \lambda \frac{\partial g}{\partial \b^{(j)}_i} + \lambda_j
\end{equation}
to eliminate the hyperplane multipliers $\lambda_j$ from the system (\ref{system}).

Using ths procedure in the codimension one case, the number of non-constant rows in $M_{\lambda}$ is equal to the number of additional top-dimensional simplices created in the subdivision.
\end{remark}

\begin{remark}\label{alwayslift}
If the top-dimensional homology of the simplicial complex is trivial, then  $\ker \partial=0$ so that the ellipsoid $g=0$ has codimension one.  Hence, every root of $Q(\lambda)$ lifts to a solution of the Lagrange multiplier system (\ref{system}).  In this case the maximum root is the maximum value of the functional $f=|Lz|^2$.
\end{remark}

\begin{figure}
\includegraphics[scale=.75, bb=176.25375 346.8163 354.92375 518.22375]{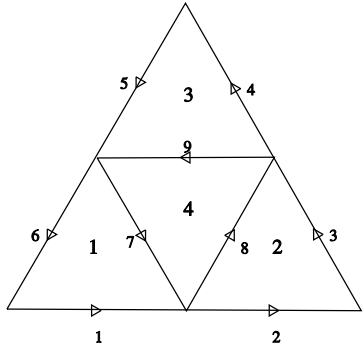}
\caption{Non-stellar subdivision of the 2-simplex}\label{nonstellarex}
\end{figure}

\begin{ex}\label{extremalex}
Consider the non-stellar subdivision in Figure \ref{nonstellarex} where the 2-simplices are given the standard counter-clockwise orientation.  Let $z=(b_i)$ be the coordinates for the space of $2$-chains relative the 2-simplices.  The ellipsoid and hyperplane conditions are given by
\begin{align}
g &= |\partial z|^2-1 = 3(b_1^2+b_2^2+b_3^2+b_4^2)-2b_4(b_1+b_2+b_3)-1 =0\\
h &= b_1+b_2+b_3+b_4 = 0
\end{align}
respectively, and the functional is
\begin{equation}
f=|Lz|^2 = (b_1+b_2)^2+(b_1+b_3)^2+(b_2+b_3)^2
\end{equation}
The Lagrange multiplier condition $\nabla f-\lambda\nabla g=\lambda_1 \nabla h$ yields four more equations:
\begin{align}
2(2b_1+b_2+b_3)-2\lambda(3b_1-b_4) &= \lambda_1\\
2(b_1+2b_2+b_3)-2\lambda(3b_2-b_4) &= \lambda_1\\
2(b_1+b_2+2b_3)-2\lambda(3b_3-b_4) &= \lambda_1\\
-2\lambda(3b_4-b_1-b_2-b_3) &= \lambda_1
\end{align}
As noted in Remark \ref{alternative}, we may use any of these to eliminate $\lambda_1$ from the remaining three to get a linear system over $b_i$ with coefficients at most linear in $\lambda$.  For example, using the last equation to eliminate $\lambda_1$ and combining with the hyperplane condition, we have
\begin{equation}
M_{\lambda}z=
\begin{pmatrix}
4-8\lambda & 2-2\lambda & 2-2\lambda & 8\lambda\\
2-2\lambda & 4-8\lambda & 2-2\lambda & 8\lambda\\
2-2\lambda & 2-2\lambda & 4-8\lambda & 8\lambda\\
1 & 1 & 1 & 1
\end{pmatrix}
\begin{pmatrix}
b_1\\ b_2\\ b_3\\ b_4
\end{pmatrix}
=0
\end{equation}
where $z=(b_i)$. Then,
\begin{equation}
\det(M_{\lambda})=16\left(\lambda-\frac{2}{9}\right)\left(\lambda-\frac{1}{3}\right)^2
\end{equation}
and by Remark \ref{alwayslift}, the desired extremum in codimension one is $1/3$.

For higher codimension, we need to compute additional hyperplane conditions.  Now, $V_2$ is by definition normal to the hyperplane $h$ representing $\im i_*$.  Hence, it is the three-dimensional space orthogonal to the line spanned by $(1,1,1,1)$, and we may arbitrarily choose a basis, say
\begin{equation*}
v_1=(1,0,0,-1)\quad\quad v_2=(0,1,0,-1)\quad\quad v_3=(0,0,1,-1)
\end{equation*}
relative the basis of 2-simplices.  Using the oriented edges in the figure as a basis for the $1$-chain space, we find that
\begin{align*}
\partial v_1 &= (1,0,0,0,0,1,-2,-1,-1)\\
\partial v_2 &= (0,1,1,0,0,0,-1,-2,-1)\\
\partial v_3 &= (0,0,0,1,1,0,-1,-1,-2)
\end{align*}
are linearly independent so that the vectors $\gamma_l=\partial v_l$ are a basis for $\partial V_2$.  If we define hyperplanes $Q_l$ to be perpendicular to $\gamma_l$, we may set $\nabla Q_l=\gamma_l$.  Combining with the party relations $P_j$ in the $1$-chain space, the Lagrange system for extremizing in codimension two is
\begin{equation}
\nabla f-\lambda\nabla g=A\mu
\end{equation}
where $\mu=(\lambda_1,\ldots,\lambda_6)^T$ and
\begin{equation}
A=\begin{pmatrix}\nabla P_j & \nabla Q_l\end{pmatrix}=
\begin{pmatrix}
1 & 0 & 0 & 1 & 0 & 0\\
1 & 0 & 0 & 0 & 1 & 0\\
0 & 1 & 0 & 0 & 1 & 0\\
0 & 1 & 0 & 0 & 0 & 1\\
0 & 0 & 1 & 0 & 0 & 1\\
0 & 0 & 1 & 1 & 0 & 0\\
0 & 0 & 0 & -2 & -1 & -1\\
0 & 0 & 0 & -1 & -2 & -1\\
0 & 0 & 0 & -1 & -1 & -2
\end{pmatrix}
\end{equation}
Since the columns are independent, $A$ has a left-inverse, and we may eliminate the variables $\lambda_i$ and proceed with the extremization algorithm as in the codimension one case.
\end{ex}

\section{Elementary Stellar Subdivisions}\label{stellar}

In this section, we use the geometric interpretation of Subsection \ref{geometric} to obtain bounds for the extremal value in codimension one for elementary stellar subdivision of an isolated $d$-simplex or along an interior $k$-simplex of a $d$-complex.  We proceed as follows: First, we show that for an elementary stellar subdivision, the functional $|Lz|$ reduces to the norm $|z|$ over the intersection locus, so that we are in effect extremizing the norm (Theorem \ref{norm}).  Second, we choose well-adapted bases for studying the spectrum of the symmetric matrix defining the (conic) intersection locus.  The reciprocal of the smallest eigenvalue maximizes the square-norm, $|z|^2$, and is the square of our desired invariant.

\begin{prop}\label{norm}
For an elementary stellar subdivision of a simplicial complex, $|Lz|=|z|$ on the intersection locus.  Thus, the functional extremized is precisely the norm on the intersection locus.
\end{prop}
\begin{proof}
We need only work in the neighborhood of simplices incident to the subdivided simplex.  Hence, all references in the argument will be with respect to the star neighborhood of the subdivided simplex.  Let $z=\sum b_i\b_i$ where $\b_i$ are the top-dimensional simplices.  Note that each codimension one party is incident to every member of incident top-dimensional parties.  Thus, in the product $Lz$, the rows corresponding to interior codimension one parties are of the form $\sum \pm P_j$ for incident parties $P_j$.  Restricting to the intersection locus, these vanish by the hyperplane condition.  Boundary codimension one parties (that is, those in the boundary of the neighborhood) are incident only to one top-dimensional member in a bijective manner so that $|Lz|=|z|$.
\end{proof}

\begin{figure}
\includegraphics[scale=.75, bb=176.25375 346.8163 354.92375 518.22375]{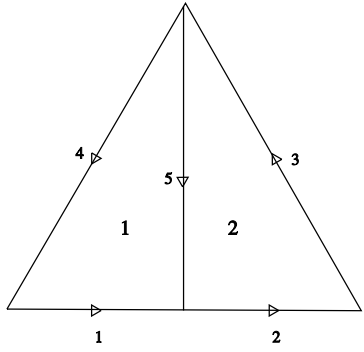}
\caption{Stellar subdivision of the 2-simplex}\label{stellarex}
\end{figure}

\begin{ex}
Consider the elementary stellar subdivision of the $2$-simplex in Figure \ref{stellarex}.  In this example we will compute relative an orientation different from that most natural.  Of course, the result is independent of such choices of orientation.

Let $z=(b_i)$ be coordinates for the $2$-chains relative the faces $f_i$ with $f_1$ oriented clockwise and $f_2$ counter-clockwise.  Note that under this orientation, the sole $2$-party is described by $f_2-f_1$ so that $b_2-b_1=0$.  If $e_i$ are the oriented edges in the diagram, we set $\{e_1+e_2,e_3,e_4\}$ as an ordered basis for $1$-chains. With respect to these two bases, the map $L$ is
\begin{equation*}
L\begin{pmatrix}b_1\\b_2\end{pmatrix}
=\begin{pmatrix}-1 & 1\\0 & 1\\1 & 0\end{pmatrix}
\begin{pmatrix}b_1\\b_2\end{pmatrix}
=\begin{pmatrix}b_2-b_1\\b_2\\b_1\end{pmatrix}
=\begin{pmatrix}0\\b_2\\b_1\end{pmatrix}
\end{equation*}
so that $|Lz|=|z|$ as claimed.
\end{ex}

For convenience, we henceforth assume that the subdivided complex is oriented so that for each party, the boundary of the sum of its members is the boundary of the party.  We call this the \emph{induced orientation} of the subdivision.  Now we change basis so that the hyperplane conditions correspond to certain axial-directions.  Since hyperplanes are of the form $\sum \b_i$, for cetain $\b_i$, we complete to a basis starting with vectors of this form.  Our convenient choice of basis is an orthogonal set when restricted to the subspace perpendicular to $\im i_*$, and this facilitates the required eigenvalue computations.

\subsection{Case of an Isolated $d$-Simplex}
The following matrices are helpful in our analysis:

\begin{defn}
Let $k>0$ and $m$ any real number.  Define $k\times k$ matrices $J_k$ and $A_{m,k}$ as follows
\begin{equation}
J_k=
\begin{pmatrix}
1 & -1 & \cdots & -1\\
\vdots & & I_{k-1} &\\
1 & & &
\end{pmatrix}
\quad\quad
A_{m,k}=
\begin{pmatrix}
m & & -1\\
& \ddots & \\
-1 &  & m
\end{pmatrix}
=(m+1)I_k-\mathbbm{1}
\end{equation}
\end{defn}

We sometimes omit the subscript for $J_k$ when the dimension is clear.  Next, by a $(k,l)$-minor of a matrix $A$, we mean the submatrix obtained by omitting the $k$-th row and $l$-th column.

\begin{prop}\label{changebasis}
The columns of $J_k$ constitute a basis for $\R^k$, and $J^{-1}A_{m,k}J$ is diagonal.  Furthermore, the $(1,1)$-minor of the latter matrix is $(m+1)I_{k-1}$.  In other words, the endomorphism $A_{m,k}$ is decomposes over $\langle \mathbbm{1}\rangle\oplus \langle \mathbbm{1}\rangle^{\perp}$ and restricts to a homothety on the latter summand. 
\end{prop}
\begin{proof}
A direct computation shows that
\begin{equation}
J^{-1}= \frac{1}{k}
\begin{pmatrix}
1 & 1 & \cdots & 1\\
-1 & k-1 & & -1\\
\vdots & & \ddots & \\
-1 & -1&  & k-1
\end{pmatrix}
\end{equation}
where the lower $(k-1)\times(k-1)$-minor is $kI_{k-1}-\mathbbm{1}$.  Alternatively, inspection shows that the first column $\mathbbm{1}$ is orthogonal to the remaining which are mutually perpendicular.  The columns of $J$ are clearly eigenvectors of $A_{m,k}$, and the eigenvalues on $\langle \mathbbm{1}\rangle^{\perp}$ are easily seen to be identical.
\end{proof}

\begin{cor}\label{sameeigenvalues}
Let $Q$ be an orthogonal matrix whose first column is parallel to $\mathbbm{1}$.  Then, the $(1,1)$-minors of $Q^TA_{m,k}Q$ and $J^{-1}A_{m,k}J$ have the same eigenvalues.
\end{cor}
\begin{proof}
The first columns of $J$ and $Q$ are parallel, and the remaining columns of $J$ are orthogonal to $\langle \mathbbm{1}\rangle$.  Thus, both $(1,1)$-minors are representing the same endomorphism on $\langle \mathbbm{1}\rangle^{\perp}$.
\end{proof}

\begin{ex}
For the case $k=4=m$, we have
\begin{equation*}
J_4=
\begin{pmatrix}
1 & -1 & -1 & -1\\
1 & 1 & 0 & 0\\
1 & 0 & 1 & 0\\
1 & 0 & 0 & 1
\end{pmatrix}
\quad\quad
A_{4,4}=
\begin{pmatrix}
4 & -1 & -1 & -1\\
-1 & 4 & -1 & -1\\
-1 & -1 & 4 & -1\\
-1 & -1 & -1 & 4
\end{pmatrix}
\end{equation*}
with
\begin{equation*}
J_4^{-1}A_{4,4}J_4=
\begin{pmatrix}
1 & 0 & 0 & 0\\
0 & 5 & 0 & 0\\
0 & 0 & 5 & 0\\
0 & 0 & 0 & 5
\end{pmatrix}
\end{equation*}
so that omitting the first row and column yields $(4+1)I_{4-1}$.
\end{ex}

We now relate the above constructions to the elementary stellar subdivision of an isolated $d$-simplex.

\begin{lem}\label{isolatedN}
For an elementary stellar subdivision of an isolated $d$-simplex along a $k$-simplex, we have $\partial^*_{top}\partial_{top}=A_{d+1,k+1}$.
\end{lem}
\begin{proof}
In the subdivided complex, each top-dimensional simplex has $(d+1)$-faces of codimension one and is incident to $k$ top-dimensional simplices.  The result follows by definition of $\partial^*$.
\end{proof}

\begin{prop}[Isolated Simplex]\label{isolatedcase}
For an elementary stellar subdivision of an isolated $d$-simplex along a $k$-simplex, the codimension one invariant value is
\begin{equation*}
C_{d-1} = \frac{1}{\sqrt{d+2}}
\end{equation*}
In particular, it is identical for all stellar subdivisions of an isolated $d$-simplex.
\end{prop}
\begin{proof}
By Proposition \ref{norm}, we are seeking the maximum norm on the intersection of the ellipsoid $x^T\partial^*\partial x=1$ and the hyperplane determined by $\im i_*$.  Using the induced (consistent) orientation on the top-simplices, the latter has perpendicular $\mathbbm{1}$.  We change to a new basis by using any orthogonal matrix $Q$ whose first column is parallel to $\mathbbm{1}$.  Then, the intersection locus is given by the vanishing of the first coordinate with respect to the new basis, namely $y^TMy=1$ where $M$ is the $(1,1)$-minor of $Q^T\partial^*\partial Q$.  By Lemma \ref{isolatedN}, $\partial^*\partial=A_{d+1,k+1}$, so that by Proposition \ref{changebasis} and its corollary (\ref{sameeigenvalues}), the eigenvalues of $M$ are those of $(d+2)I_k$.  Hence, the maximum norm on the intersection locus is $1/\sqrt{d+2}$.
\end{proof}

\begin{remark}\label{insufficient}
The extremal values do not differentiate every subdivision.  Combining the preceding proposition with Example \ref{C0stellarcase}, we see that the sequence of invariants $C_k$ cannot distinguish between the two elementary stellar subdivisions of an isolated $2$-simplex.
\end{remark}

\subsection{Case of Subdividing along an Interior $k$-Simplex}
Now, we consider a stellar subdivision of a $d$-complex along an interior $k$-simplex.  We will order top-dimensional simplices in within each party and then collect together.  Then, with respect to this adapted ordering, the matrix for $\partial^*_{top}\partial_{top}$ consists of square $(k+1)$-blocks of either $0,-I,A_{d+1,k+1}$.  The diagonal blocks are $A_{d+1,k+1}$ and reflect incident relations within a party; the $-I$ blocks reflect incidence relations between parties, since for $k<d$, if two top-dimensional parties are incident then there is an incidence bijection between the respective party members.  Instead of $J$, we use a change of basis matrix $P$ that is block diagonal with blocks $J_{k+1}$.  Thus, the hyperplane relations are encoded by columns numbered $1+j(k+1)$ for $j$ a non-negative integer, and the intersection locus is described by the minor obtained by omitting these columns and the same-numbered rows.  By the definition of block multiplication, the consequent minor of the matrix $P^{-1}\partial^*\partial P$ then consists of square $k$-blocks of $0,-I,(d+2)I$ where the $(d+2)I_k$ appear along the diagonal and the $-I_k$ blocks are in the same relative positions as $-I_{k+1}$ in $\partial^*\partial$.  Using an argument analogous to that in Proposition \ref{isolatedcase}, the eigenvalues for the minor of $P^{-1}\partial^*\partial P$ are the same as that for $Q^T\partial^*\partial Q$ where $Q$ is block-diagonal with orthogonal blocks having their first columns each parallel to $\mathbbm{1}_{k+1}$.

At this point, Gersgorin's Theorem tells us that a lower bound for the eigenvalues is $(d+2)-(d-k)$ where $(d-k)$ is the number of entries of $-1$ in each row.  Recall that these $-1$ entries reflect incidence relations between non-singular parties.  Thus, an intuitive reason for this number is that at the added vertex, $k$ of the $d$ directions are occupied by the subdivided $k$-simplex, leaving only $(d-k)$ directions which are all filled by top-dimensional neighbors because $v$ is an interior point.  This lower bound of $(k+2)$ is indeed assumed by the eigenvector $\mathbbm{1}$.  Thus, we may informally conclude that $C_{d-1}=1/\sqrt{k+2}$, for an elementary stellar subdivision on an interior $k$-simplex.

More formally, we may proceed as follows:

\begin{defn}
Let $j,k,l$ be non-negative integers.  Given a square matrix $A$ of dimension $k$, we define a sequence of square matrices $R^j(A)$ of size $k2^j$ for $j\geq 0$ as follows:
\begin{equation}
R^0(A)= A \quad\quad R^{j+1}(A) = \begin{pmatrix}R^j & -I\\ -I & R^j\end{pmatrix}
\end{equation}
For $l,d\geq 0$ and $k\leq d$, we will write
\begin{equation}
M_{d,k}^{(l)} = R^l((d+2)I_k) \quad\quad N_{d,k}^{(l)} = R^l(A_{d+1,k+1})
\end{equation}
\end{defn}

\begin{ex}
Setting 
\begin{equation*}
R^0(A) = A = \begin{pmatrix}a & b\\c &d\end{pmatrix}
\end{equation*}
we have
\begin{equation*}
R^1(A)= \begin{pmatrix}a & b & -1 &\\c & d &  & -1\\-1 &  & a & b\\ & -1 & c & d\end{pmatrix}
\end{equation*}
\end{ex}

The positive integer $l$ in $R^l(A)$ thus counts the number of ``binary levels'' in the resulting matrix.

Thus, Proposition \ref{changebasis} may be restated as follows: the $(1,1)$-minor of $J^{-1}N_{d,k}^{(0)}J$ is $M_{d,k}^{(0)}$.  Also, the block-multiplication exhibited in the example
\begin{equation*}
\begin{pmatrix}J^{-1} & \\ & J^{-1}\end{pmatrix}
\begin{pmatrix}N_{d,k}^{(0)} & -I\\ -I & N_{d,k}^{(0)}\end{pmatrix}
\begin{pmatrix}J & \\ & J\end{pmatrix}
=\begin{pmatrix}J^{-1}N_{d,k}^{(0)}J & -I\\ -I & J^{-1}N_{d,k}^{(0)}J\end{pmatrix}
\end{equation*}
shows that more generally, if $P$ is block-diagonal with blocks $J$, then deleting rows and columns in $P^{-1}N_{d,k}^{(l)}P$ numbered $1+j(k+1)$ for $j$ a non-negative integer,  yields $M_{d,k}^{(l)}$.

\begin{prop}[Interior Case]\label{interiorcase}
For an elementary stellar subdivision of a $k$-simplex in the interior of a $d$-complex, we have $\partial^*_{top}\partial_{top}=N_{d,k}^{(d-k)}$, and the square-norm of points on the associated intersection locus is sharply bounded by
\begin{equation}
\frac{1}{2d+2-k}\leq |z|^2 \leq \frac{1}{k+2}
\end{equation}
In particular, for such subdivisions,
\begin{equation*}
C_{d-1} = \frac{1}{\sqrt{k+2}}
\end{equation*}
\end{prop}
\begin{proof}
Recall that by the local nature of an elementary stellar subdivision, we only need to work with top-simplices incident to the subdivided $k$-simplex.  First, we order the top-simplices within each party and then collect as follows: Choose an ordering of the $d$-coordinate directions.  Pick a top-party whose members begin the order.  Next, collect the top-parties incident to this and order their members.  Then, collect the top-parties incident to those already numbered and order them so that they are in the same relative order as the parties to which they are incident.  Continue likewise until all $d$-directions are used.  Each step adds twice the number of parties to the running total.  For the change of basis matrix, we use the $(k+1)2^{d-k}$-dimensional block diagonal matrix
\begin{equation*}
P=\diag(J_{k+1},\ldots,J_{k+1})
\end{equation*}
so that the first column of each block reflects a party hyperplane relation.  Next,
\begin{equation}
\partial_{top}^T\partial_{top}=N_{d,k}^{(d-k)}
\end{equation}
on the neighborhood of top-simplices incident to the subdivided interior $k$-simplex.  To see this note that the diagonal blocks of $A_{d+1,k+1}$ in $N_{d,k}^{(d-k)}$ reflect incidence relations in each party, and successive applications of the operator $R$ pick up the incident parties in each of the remaining $(d-k)$ directions for which a top-simplex may be incident another in the relevant neighborhood.

Letting $Q$ be a block-diagonal matrix of orthogonal blocks with first column parallel to $\mathbbm{1}_{k+1}$, the intersection locus is thus given by a minor of $Q^T N_{d,k}^{(d-k)} Q$ where rows and columns numbered $1+j(k+1)$ for $j$ a non-negative integer, are omitted.  By an argument analogous to that in Corollary \ref{sameeigenvalues}, we see that the desired eigenvalues are identical to those of the corresponding minor of $P^{-1}N^{(d-k)}_{d,k}P$, namely $M_{d,k}^{(d-k)}$.  

Next, bounds for the square-norm on the conic intersection locus are given by the reciprocals of the extreme values for the eigenvalues of $M_{d,k}^{(d-k)}$.  Since there are $(d-k)$ entries of $-1$ in each row, by Gersgorin's Theorem, the eigenvalues $\lambda$ are bounded by
\begin{equation*}
k+2 \leq \lambda \leq 2d+2-k
\end{equation*}
The lower bound is clearly achieved by the eigenvector $\mathbbm{1}=(1,\ldots,1)$.  To show that the upper bound is also attained, define $v_{k,0}=(1,\ldots,1)$ a $k$-vector, and $v_{k,i+1}=v_{k,i}\oplus -v_{k,i}$, a $k2^{i+1}$-vector formed by ``concatenation'' of two $k2^i$-vectors $v_{k,i}$.  For example, $v_{k,1}=(1,\ldots,1,-1,\ldots,-1)$ a $2k$-vector.  Note that
\begin{align*}
M_{d,d-1}^{(1)}v_{d-1,1} &=
\begin{pmatrix}(d+2)I & -I\\-I & (d+2)I\end{pmatrix}
\begin{pmatrix}\mathbbm{1}\\-\mathbbm{1}\end{pmatrix}
= (2d+2-(d-1))v_{d-1,1}\\
M_{d,d-2}^{(2)}v_{d-1,2} &=
\begin{pmatrix}(d+2)I & -I & -I &\\-I & (d+2)I &  & -I\\-I & & (d+2)I & -I\\ & -I & -I & (d+2)I\end{pmatrix}
\begin{pmatrix}\mathbbm{1}\\-\mathbbm{1}\\ -\mathbbm{1}\\ \mathbbm{1}\end{pmatrix}\\
&= (2d+2-(d-2))v_{d-1,2}\\
\end{align*}
and proceeding iteratively by induction on depth, the upper bound is attained.  Finally, interpreting geometrically, the bound on the square-norm follows.
\end{proof}

\begin{ex}
Consider an elementary stellar subdivision of an interior standard $3$-simplex $\Delta_3$ on an interior $1$-simplex.  Grouping the top-dimensional simplices by parties, we have
\begin{equation*}
\partial_{top}^*\partial_{top}=N_{3,1}^{(2)}=
\begin{pmatrix}
4 & -1 & -1 & 0 & -1 & 0 & 0 & 0\\
-1 & 4 & 0 & -1 & 0 & -1 & 0 & 0\\
-1 & 0 & 4 & -1 & 0 & 0 & -1 & 0\\
0 & -1 & -1 & 4 & 0 & 0 & 0 & -1\\
-1 & 0 & 0 & 0 & 4 & -1 & -1 & 0\\
0 & -1 & 0 & 0 & -1 & 4 & 0 & -1\\
0 & 0 & -1 & 0 & -1 & 0 & 4 & -1\\
0 & 0 & 0 & -1 & 0 & -1 & -1 & 4
\end{pmatrix}
\end{equation*}
with change of basis matrix a block-diagonal matrix $P$ with diagonal blocks:
\begin{equation*}
J=\begin{pmatrix} 1 & -1\\ 1 & 1\end{pmatrix}
\end{equation*}
Then,
\begin{equation*}
P^{-1}NP =
\begin{pmatrix}
3 & 0 & -1 & 0 & -1 & 0 & 0 & 0\\
0 & 5 & 0 & -1 & 0 & -1 & 0 & 0\\
-1 & 0 & 3 & 0 & 0 & 0 & -1 & 0\\
0 & -1 & 0 & 5 & 0 & 0 & 0 & -1\\
-1 & 0 & 0 & 0 & 3 & 0 & -1 & 0\\
0 & -1 & 0 & 0 & 0 & 5 & 0 & -1\\
0 & 0 & -1 & 0 & -1 & 0 & 3 & 0\\
0 & 0 & 0 & -1 & 0 & -1 & 0 & 5
\end{pmatrix}
\end{equation*}
and the minor describing the intersection locus is
\begin{equation*}
M_{3,1}^{(2)} =
\begin{pmatrix}
5 & -1 & -1 & 0\\
-1 & 5 & 0 & -1\\
-1 & 0 & 5 & -1\\
0 & -1 & -1 & 5
\end{pmatrix}
\end{equation*}
Furthermore, we see that 
\begin{equation*}
M_{3,1}^{(2)} \mathbbm{1} =
\begin{pmatrix}
5 & -1 & -1 & 0\\
-1 & 5 & 0 & -1\\
-1 & 0 & 5 & -1\\
0 & -1 & -1 & 5
\end{pmatrix}
\begin{pmatrix}
 1\\ 1\\ 1\\ 1
\end{pmatrix}
=3\begin{pmatrix} 1\\ 1\\ 1\\ 1\end{pmatrix}
\end{equation*}
and
\begin{equation*}
M_{3,1}^{(2)} v_{1,2} =
\begin{pmatrix}
5 & -1 & -1 & 0\\
-1 & 5 & 0 & -1\\
-1 & 0 & 5 & -1\\
0 & -1 & -1 & 5
\end{pmatrix}
\begin{pmatrix}
 1\\ -1\\ -1\\ 1
\end{pmatrix}
=7\begin{pmatrix} 1\\ -1\\ -1\\ 1\end{pmatrix}
\end{equation*}
where $v_{1,2}=(1,-1,-1,1)$.
\end{ex}

\begin{remark}\label{isolatedstellar}
When the subdivided $k$-simplex is on the boundary of the $d$-complex, the analysis is similar in spirit to that above.  In the codimension one case, the local nature of the elementary stellar subdivision means that this case is equivalent to an elementary stellar subdivision of an isolated $d$-simplex treated above.  However, if the codimension of the subdivided simplex is greater than one, the intersection locus is given by a minor of $M_{d,k}^{(d-k)}$ since some of the incidence relations are not present.  By Gersgorin's Theorem, the bounds on the eigenvalues remain valid albeit not necessary sharp; namely, if $\mu$ denotes the maximum number of incident top-simplices in the relevant neighborhood, then
\begin{equation*}
d+2-\mu \leq \lambda\leq d+2+\mu
\end{equation*}
with associated bounds on $C_{d-1}$.
\end{remark}

The computed values for $C_{d-1}$ in Propositions \ref{isolatedcase} and \ref{interiorcase} are consistent with the combinatorial bound of Proposition \ref{lowerbound}.  In the case of an isolated $d$-simplex, we have $F_{\s}=2$ and $N_{\s}\leq 2$ so that as in Example \ref{combmfdcase}, we have
\begin{equation*}
C_{d-1}\geq \frac{1}{\sqrt{d+2}}
\end{equation*}
actually a sharp bound.  For the case of subdividing an interior $k$-simplex, $F_{\s}=2$ and $N_{\s}=2$.  To see the latter equality, note that because of the hypothesis on $\s$, it is incident to exactly two $d$-simplices.  Each of these $d$-simplices has an unique $(d-1)$-face that coincides with a party.  The other codimension one faces support the added vertex and hence if supported by a party, that party is not singly-represented.  Thus, we obtain
\begin{equation*}
C_{d-1}\geq \frac{1}{\sqrt{d+3}}
\end{equation*}
and the inequalities $k+2\leq d+2\leq d+3$ assure consistency.

\begin{remark}
In the case of an isolated simplex, $C_{d-1}$ reflects only global information, whereas in the case of subdividing an interior simplex, $C_{d-1}$ reflects only local information.
\end{remark}

\section{Appendix: The Laplacian on Chain Complexes}\label{laplheat}

In this appendix we provide some basic properties of the laplacian of a chain complex.

Let $\calC_*=\{(V_*,g),\partial_*\}$ be a chain complex of finite-dimensional inner product spaces:
\begin{equation}
\cdots \rightleftharpoons V_{k+1} \rightleftharpoons V_k \stackrel{\partial_k}{\rightleftharpoons} V_{k-1} \rightleftharpoons \cdots
\end{equation}
with adjoints $\partial^*$.
\begin{defn}
The \emph{Laplacian} for a chain complex $\calC_*=\{(V_*,g),\partial_*\}$ of finite-dimensional inner product spaces is 
\begin{equation}
\laplacian_k = \partial_k^*\partial_k + \partial_{k+1}\partial_{k+1}^*
\end{equation}
\end{defn}

\begin{prop}\label{basic}
The Laplacian is non-negative, self-adjoint, and satisfies
\begin{equation}
\dim \ker \laplacian_k = \dim H_k(\calC_*)
\end{equation}
\end{prop}
\begin{proof}
Self-adjointness is evident from the definition.  For any $v\in V_k$, we have
\begin{equation}\label{laplacian-relation}
(\laplacian_k v,v) = |\partial_k v|^2 + |\partial_{k+1}^* v|^2 \geq 0
\end{equation}
where the norms are those relative the respective inner products.  This shows that all eigenvalues are non-negative.

Finally, to show that $\dim \ker \laplacian_k = \dim H_k(\calC_*)$, note that (\ref{laplacian-relation}) implies that $v\in \ker \laplacian_k$ iff $v\in \ker \partial_k\cap \ker\partial^*_{k+1}$.  Now by the definition of adjoint, $\ker \partial_{k+1}^* = (\im \partial_{k+1})^{\perp}$ relative the inner product on $V_{k+1}$.  Thus,
\begin{equation*}
\ker \laplacian_k = \ker \partial_k \cap (\im\partial_{k+1})^{\perp}
\end{equation*}
and since $\im \partial_{k+1}\subset \ker\partial_k$, we have
\begin{equation*}
\dim\ker\laplacian_k = \dim\left( \ker \partial_k \cap (\im\partial_{k+1})^{\perp} \right) = \dim(\ker\partial_k) - \dim(\im \partial_{k+1})
\end{equation*}
as desired.
\end{proof}

By standard linear algebra, there is a \emph{Hodge Decomposition} (see \cite{eckmann00}):
\begin{equation}
V_k = \ker \laplacian_k \oplus \im \partial_{k+1} \oplus \im \partial_k^*
\end{equation}
which is much easier to prove in this setting than that of differential forms because the inner product spaces are assumed finite-dimensional.

\bibliographystyle{plain}
\bibliography{subdivision}

\end{document}